
\input amssym.def
\input amssym.tex

\font\eufrak=eufm10

\def\D{\mathop{\rm{}D}\nolimits}
\def\F{{\Bbb F}}
\def\M{\mathop{\rm{}M}\nolimits}
\def\T{\mathop{\rm{}T}\nolimits}
\def\R{{\Bbb R}}

\def\Sym{\mathop{\rm{}S}\nolimits}
\def\SL{\mathop{\rm SL}\nolimits}
\def\SD{\mathop{\rm SD}\nolimits}
\def\GL{\mathop{\rm GL}\nolimits}
\def\diag{\mathop{\rm diag}\nolimits}
\def\th{\mathop{\rm th}\nolimits}
\def\mod{\mathop{\rm mod}\nolimits}

\magnification=\magstep1
\tolerance=300
\pretolerance=200
\hfuzz=1pt
\vfuzz=1pt
\parindent=35pt
\mathsurround=1pt
\parskip=1pt plus .25pt minus .25pt
\baselineskip=11pt
\normallineskiplimit=1pt
\abovedisplayskip=10pt plus 2.5pt minus 7.5pt
\abovedisplayshortskip=0pt plus 2.5pt
\newdimen\jot
\jot=2.5pt

\hoffset=0.3in
\voffset=-1\baselineskip
\hsize=5.8 true in
\vsize=53.25\baselineskip

\font\eightrm=cmr8

\font\bfone=cmbx10 scaled\magstep1
\font\smc=cmcsc10

%
%
\mathchardef\emptyset="001F

\def\nin{\noindent}

\nopagenumbers
\def\rightheadline{\hfil\smc\lastname\hfil\tenbf\folio}
\def\leftheadline{\tenbf\folio\hfil\smc\lastname\hfil}
\def\firstheadline{\vbox{\baselineskip=8pt
    \parindent 0pt \obeylines \eightrm
    Semigroup Forum  \hfill{\date}
    \par\nin Springer-Verlag New York Inc.}}
\headline={\ifnum\pageno=1\firstheadline
           \else\ifodd\pageno\rightheadline
           \else\leftheadline\fi\fi}
\def\datum{\ifcase\month\or January\or February\or March\or April
    \or May\or June\or July\or August\or September\or October
    \or November\or December\fi\space\number\day, \number\year}
\def\date{{\eightrm Version of~}\datum}

\def\title{AUTHOR: INSERT TITLE!!}
\def\author{AUTHOR: INSERT YOUR NAME(S)}
\def\lastname{AUTHOR: INSERT YOUR LASTNAME(S) (WITHOUT INTIALS)}
\def\editor{AUTHOR: INSERT EDITOR'S NAME}
\def\thanks#1{\footnote*{\eightrm#1}}

\def\title{\centerline{\bfone\titleone}\centerline{\bfone\titletwo}}
\def\titleone{}
\def\titletwo{}
\def\firstpage{\vbox{\vskip1.5truecm}
       \title
       \vskip.75truecm\centerline{\bf\author}
       \medskip
       \centerline{\eightrm Communicated by \editor}
       \vskip.75truecm\rm}

\def\sectionheadline #1{\vskip-\lastskip\bigbreak
                      \centerline{\bf #1}\nobreak\medskip\nobreak}

        \newtoks\literat
\def\[#1 #2\par{\literat={#2\unskip.}%
        \hbox{\vtop{\hsize=.1\hsize\nin [#1]\hfill}
        \vtop{\hsize=.9\hsize\nin\the\literat}}\par
        \vskip.3\baselineskip}

\def\references{\sectionheadline{\bf References}\frenchspacing
                \entries\par}

\def\address{Author: {\tt$\backslash$def$\backslash$address$\{$??$\}$}}
\def\addresstwo{}
\def\lastpage{\par\vbox{\vskip.8truecm\nin
\line{\vtop{\hsize=.45\hsize{\eightrm\parindent=0pt\baselineskip=8pt
            \nin\address}}\qquad
      \vtop{\hsize=.45\hsize\nin{\eightrm\parindent=0pt
            \baselineskip=8pt\addresstwo}}\hfill}
}}

\def\qed{{\unskip\nobreak\hfil\penalty50\hskip .001pt \hbox{}\nobreak\hfil
          \vrule height 1.2ex width 1.1ex depth -.1ex
          \parfillskip=0pt\finalhyphendemerits=0\medbreak}\rm}


\def\Proposition #1. {\bigbreak\vskip-\parskip\noindent{\bf Proposition #1.}
    \quad\it}
\def\Propositions. {\bigbreak\vskip-\parskip\noindent{\bf Proposition.}
    \quad\it}
\def\Theorem #1. {\bigbreak\vskip-\parskip\noindent{\bf  Theorem #1.}
    \quad\it}
\def\Corollary #1. {\bigbreak\vskip-\parskip\nin{\bf Corollary #1.}
    \quad\it}
\def\Lemma #1. {\bigbreak\vskip-\parskip\noindent{\bf  Lemma #1.}\quad\it}
\def\Lemmad. {\bigbreak\vskip-\parskip\noindent{\bf  Lemma.}\quad\it}
\def\Definition #1. {\rm\bigbreak\vskip-\parskip\noindent{\bf Definition #1.}
    \quad}
\def\Remark #1. {\rm\bigbreak\vskip-\parskip\noindent{\bf Remark #1.}\quad}
\def\Remarkd. {\rm\bigbreak\vskip-\parskip\noindent{\bf Remark.}\quad}
\def\Exercise #1. {\rm\bigbreak\vskip-\parskip\noindent{\bf Exercise #1.}
    \quad}
\def\Example #1. {\rm\bigbreak\vskip-\parskip\noindent{\bf Example #1.}\quad}
\def\Examples #1. {\rm\bigbreak\vskip-\parskip\noindent{\bf Examples #1.}\quad}
\def\Proof#1.{\rm\par\ifdim\lastskip<\bigskipamount\removelastskip\fi
    \smallskip\noindent{\bf Proof.}\quad}
\def\Corollary #1. {\bigbreak\vskip-\parskip\nin{\bf Corollary.}
    \quad\it}


\hoffset=0.3in          
\voffset=-1\baselineskip 


\def\titleone{Intermediate semigroups are groups}


\def\titletwo{}


\def\author{Alexandre~A.~Panin}
 

\def\lastname{Panin}


\def\editor{Mohan Putcha}


%
%
%

\def\entries
{\[1 Borewicz, Z.I., {\it On parabolic subgroups in linear groups
over a semilocal ring,}
Vestn. Leningr. Univ., Math. {\bf 9} (1981), 187--196

\[2 Borewicz, Z.I., {\it Description of the subgroups of the full
linear group that contain the group of diagonal matrices,}
J.~Sov. Math. {\bf 17} (1981), 1718--1730

\[3 Borewicz, Z.I., and Vavilov, N.A., {\it Subgroups of the
general linear group over a semilocal ring containing the group of
diagonal matrices,}
Proc. Steklov Inst. Math. {\bf 148} (1980), 41--54

\[4 Borewicz, Z.I., and Vavilov, N.A., {\it On the definition of
a net subgroup,}
J.~Sov. Math. {\bf 30} (1985), 1810--1816

\[5 Kojbaev, V.A., {\it The subgroups of the group {\rm$\GL(2,K)$}
that contain a nonsplit maximal torus,}
J.~Math. Sci. {\bf 83} (1997), 648--653

\[6 Merzlyakov, Yu.I., ``Rational groups'',
Nauka, Moscow, 1980 (In Russian)

\[7 Pr\"ohle, P., {\it Which of the cancellative semigroups are groups?,}
Semigroup Forum {\bf 57} (1998), 438--439

\[8 Vavilov, N.A., {\it On subgroups of the general linear group
over a semilocal ring that contain the group of diagonal matrices,}
Vestn. Leningr. Univ., Math. {\bf 14} (1981), 9--15

\[9 Vavilov, N.A., {\it Subgroups of Chevalley groups containing a
maximal torus,}
Transl. Amer. Math. Soc. {\bf 155} (1993), 59--100

\[10 Vavilov, N.A., {\it Intermediate subgroups in Chevalley groups,}
Proc. Conf. Groups of Lie Type and their Geometries
{\rm(}Como -- 1993\/{\rm)}, Cambridge Univ. Press, 1995, pp.$\,$233--280}


\def\address{Dept. of Mathematics and Mechanics

St.Petersburg State University

2 Bibliotechnaya square

St.Petersburg 198904, Russia

e-mail: alex@ap2707.spb.edu

}


\def\addresstwo{

}                          


\firstpage
                                           
\sectionheadline{Abstract}

We consider the lattice of subsemigroups of the general linear
group over an Artinian ring containing the group of diagonal
matrices and show that every such semigroup is actually a group.

\sectionheadline{Introduction}

Let $R$ be an associative ring with unit, $R^*$ be its multiplicative
group, $J=J(R)$ be its Jacobson radical.
Let $G=\GL(n,R)$ be the general
linear group of degree $n$ over $R$ and $D=\D(n,R)$ be its subgroup
of diagonal matrices. For a matrix $g\in\GL(n,R)$ we denote by $g_{ij}$
its entry in the position $(i,j)$, so that
$g=(g_{ij}),\ 1\leqslant i,j\leqslant n$. As usual $g^{-1}=(g'_{ij})$
stands for the inverse of $g$, $e$~denotes the identity matrix.

Consider a square array
$\sigma=(\sigma_{ij}),\ 1\leqslant i,j\leqslant n$, consisting of
$n^2$ two-sided ideals in $R$. This array is called a {\it net
over} $R$ {\it of degree} $n$ [1], if $\sigma_{ir}\sigma_{rj}\subseteq\sigma_{ij}$
for all values of the indices $i,j,r$. A net $\sigma$ is called a
$D$--{\it net\/} if $\sigma_{ii}=R$ for all~$i$. 

For a given net $\sigma$, we put

$$M(\sigma)=\{ a=(a_{ij})\in\M(n,R): a_{ij}\in\sigma_{ij} \hbox{ for all } i,j\}.$$

The largest subgroup of
$\GL(n,R)$ contained in the multiplicative
semigroup $e+M(\sigma)$ is called the {\it net subgroup of} $\GL(n,R)$,
corresponding to the net $\sigma$~[1], and is denoted by $G(\sigma)$.

\smallskip

The lattice of subgroups of the general linear group over a field
containing the group of diagonal matrices was described by Z.I.Borewicz
in [2]. Later this result was generalised by Z.I.Borewicz and N.A.Vavilov
to the case of semilocal
rings (recall that a ring $R$ is called semilocal, if the factor-ring $R/J(R)$
is Artinian). The next theorem was proved in [8] (see [3] for a
preliminary version).

\Theorem 1.
Let $R$ be a semilocal ring such that the decomposition
of the factor-ring $R/J(R)$ in the direct sum of simple Artinian rings
does not include either fields containing less than seven elements, or
the full matrix ring $\M(2,\F_2)$. Then for every intermediate
subgroup $F,\ \D(n,R)\leqslant F\leqslant \GL(n,R)$, there exists a
unique $D$--net $\sigma$ over $R$ of degree $n$ such
that $G(\sigma)\leqslant F\leqslant N(\sigma)$, where
$N(\sigma)$ is the normaliser of the net subgroup $G(\sigma)$ in
$\GL(n,R)$.\rm

\smallskip

A number of articles (see the surveys [9], [10]) was devoted
to the description of lattices of subgroups in Chevalley groups
(or their extensions by diagonal automorphisms) containing a maximal
split torus.

\smallskip

One can try to describe the lattice of intermediate
{\it subsemigroups} of the general linear group over a ring containing
the group of diagonal matrices. In general this lattice is larger than
the lattice of intermediate {\it subgroups}, as follows from the
next result which was proved in [4].

\Theorem 2. For all nets $\sigma$ over a ring $R$
$$G(\sigma)=\GL(n,R)\cap (e+M(\sigma))$$
if and only if
for every natural number $n$ and every ideal {\eufrak a} in $R$ the
matrix ring $\M(n,R/${\eufrak a}$)$ is Dedekind finite, i.e. every
element in $\M(n,R/${\eufrak a}$)$ that has a one-sided inverse
must have a two-sided inverse.\rm

\smallskip

In the present paper we show that there is a case when
the lattice of intermediate subsemigroups is equal to the
lattice of intermediate subgroups. Hereafter an Artinian ring
stands for right or left Artinian ring.
The following result is announced in the title of the paper.

\Propositions. If $R$ is an Artinian ring, then every intermediate subsemigroup
of $\GL(n,R)$ containing $\D(n,R)$ is a subgroup.\qed

\smallskip

This proposition is a consequence of Theorem~3, which is
the main result of this paper.

If $a\in G=\GL(n,R)$, then the set of elements from the semigroup generated
by $a$ and $D=\D(n,R)$ which are of the form $\omega=\omega_1\ldots\omega_m$,
where $\omega_i=a$ for some $i$, is a semigroup denoted by
$\langle a,D\rangle$.

\Theorem 3.
Let $R$ be an Artinian ring. If $a\in G$, then $\langle a,D\rangle$
contains the identity matrix.\rm

\sectionheadline{Proof of the main theorem}

The scheme of our proof of Theorem~3 is the following.
First, we establish the result for skew fields. This part requires
some technical work. Then we immediately obtain the result for
simple and semisimple Artinian rings. At last, we complete
the proof considering ``subradical'' matrices.

\smallskip

In what follows we suppose that $n\geqslant 2$.

\smallskip

$1^o$.~Let $R=T$ be a skew field. If $G$ is a finite group, then it is obvious
that each subsemigroup of $G$ is a subgroup, so we shall assume hereafter that
$T$ is infinite.

Let $a\in\GL(n,T)$ and fix a pair of
indices $i,j$. By $a(i,j)$ we denote the matrix in $\M(n-1,T)$ obtained
from $a$ after eliminating its $i^{\th}$ row and $j^{\th}$ column.
Further, the group $\GL(n-1,T)$ is injected to the group $\GL(n,T)$ by the
correspondence
$$a\mapsto i(a)=\left(\matrix{a&0\cr
                               0&1\cr}\right)$$
We do not distinguish $a$ and its image $i(a)$, considering $a$ as an
element of $\GL(n,T)$.

\smallskip

The following assertion is of course well known (especially when $T$ is
commutative). We present its proof for the sake of completeness.

\Lemma 1.
If $a\in\GL(n,T)$, then $a'_{ji}$ is nonzero if and only if the matrix
$a(i,j)$ is invertible.
\Proof. First, let $i=j$. Without loss of generality we can assume that
$i=j=n$. There exist matrices $b_1,b_2\in\GL(n-1,T)$ such that
$b_1a(n,n)b_2$ is diagonal with the diagonal entries $d_1,\ldots,d_{n-1}$.
Let $c=b_1ab_2$ (recall our convention). Then $c'_{nn}=a'_{nn}$ and
for every $k\neq n$
$$0=\sum\limits_{r=1}^nc^{}_{kr}c'_{rn}=d^{}_kc'_{kn}+c^{}_{kn}c'_{nn}$$
It is clear that $d_k\neq 0$ for all $k$ is equivalent to $c'_{nn}\neq 0$.

The case of distinct $i$ and $j$ is easily reduced to the already
considered one, since we can multiply $a$ by the matrix of the
permutation $(ij)$.\qed

\Corollary.
If $a\in\GL(n,T)$, then there exists a permutation $\rho\in\Sym_n$
such that $a_{i\rho i}\neq 0$ for every $i$.
\Proof. Since $a$ is invertible, there exists an index $j$
such that $a^{}_{1j}a'_{j1}\neq 0$. From Lemma~1 it follows that
the matrix $a(1,j)$ is invertible. Now argue by induction on $n$.\qed

\smallskip

Now we can start to prove Theorem~3 in the case of skew fields.
Given a matrix $a\in G$, we shall find a matrix in
$\langle a,D\rangle$ which has more zero entries than $a$, and
then proceed by induction. At first we have to get rid of
``superfluous'' zeros.

\Lemma 2.
If $a\in G$, then $\langle a,D\rangle$ contains a matrix
with all diagonal entries being nonzero.
\Proof. By the Corollary from Lemma~1 there exists a permutation $\rho\in\Sym_n$
such that $a_{i\rho i}\neq 0$ for every $i$. If $d=\diag(d_1,\ldots,d_n)\in D$,
then $b_{i\rho^2i}=\sum\limits_{r=1}^na_{ir}d_ra_{r\rho^2i},$ where by $b$ we
denote the matrix $ada$. If $r=\rho i$, then the corresponding summand of this
sum is nonzero. It is clear that we can choose $d\in D$ such that the
elements $b_{i\rho^2i}$ are nonzero for all $i$.
It remains to proceed by induction on the order of $\rho$.\qed

\Lemma 3.
If $a\in G$, then there exists a matrix $b\in\langle a,D\rangle$ such that
its diagonal entries are nonzero, and if $b_{ij}=0$, then $b_{ir}b_{rj}=0$
for every $r$.
\Proof. Let's call $m^{\th}$ row of the matrix $a$ {\it good} if $a_{mm}\neq 0$
and from $a_{mj}=0$ it follows that $a_{mr}a_{rj}=0$ for every $r$.
We have to find a matrix $b\in\langle a,D\rangle$ with all rows being good.
It is clear that if $m^{\th}$ row of $a$ is good, then the same is true
for $m^{\th}$ row of each matrix from $\langle a,D\rangle$.

Suppose that $i^{\th}$ row of $a$ is not good.
It follows from Lemma~2 that we can assume $a_{ll}\neq 0$
for every $l$. Consider the set $I=\{j: a_{ij}=0\}$. We shall argue
by induction on $\vert I\vert$. By the assumption, $I$ is not
empty and $i\not\in I$. Suppose $(ada)_{ij}=0$ for every $j\in I$ and $d\in D$.
Then $a_{ir}a_{rj}=0$ for every $j\in I$ and every $r$, and we get
a contradiction ($i^{\th}$ row of $a$ is not good!). Therefore
$(ada)_{ij}\neq 0$ for some $d\in D$ and $j\in I$, and then
$a_{ir}a_{rj}\neq 0$ for some $r$. We can choose $d\in D$ such that
$(ada)_{il}\neq 0$ for all $l\not\in I$ and $l=j$, and $(ada)_{ll}\neq 0$
for every $l$. Now we can apply the induction hypothesis to get a matrix
$b\in\langle a,D\rangle$ with good $i^{\th}$ row.\qed

\Lemma 4.
If $a\in G$, then there exists a matrix $c\in\langle a,D\rangle$ such that
its diagonal entries are nonzero, and
\smallskip\indent
(i) if $c_{ij}=0$, then $c_{ir}c_{rj}=0$ for every $r$;
\smallskip\indent
(ii) if $c_{ij}\neq 0$, then $c'_{ij}\neq 0$.
\Proof. Let $b$ be a matrix whose existence is guaranteed by Lemma~3.
Consider the $D$--net $\sigma$ such that $\sigma_{ij}=0$ if and only if
$b_{ij}=0$. Then $b\in G(\sigma)$. It follows from Theorem~2 that
$b^{-1}\in G(\sigma)$. Now if a row of
$b^{-1}$ is not good, then one should
repeat the proof of Lemma~3 for $b^{-1}$ and
note that we can suppose that all diagonal matrices $d$ occurred there
have the additional property that for every nonzero entry of $b$ the
corresponding entry of $bdb$ is also nonzero.\qed

\smallskip

We are ready to prove Theorem~3 for skew fields. Without loss of generality
we can assume that the matrix $a$ satisfies the conditions $(i), (ii)$ from
Lemma~4.
Suppose that $a_{ij}=0$ for every $i<l$ and $j\neq i$. We show how to
find a matrix $b\in\langle a,D\rangle$ such that $b_{ij}=0$ for
every $i\leqslant l$ and $j\neq i$.

Let $d$ be a diagonal matrix and consider the matrix $ada$. We have
$(ada)_{ij}=\sum\limits_{r=1}^na_{ir}d_ra_{rj}$. 
Let $I=\{j: a_{lj}=0\}$. In order to obtain
$(ada)_{lj}=0$ for all $j\neq l$, we get a system of linear
equations over $T$

$$\sum\limits_{r\not\in I}\bar d_ra_{rj}=0,\ j\not\in I\cup\{l\}$$
where $\bar d_r=a_{lr}d_r,\ r\not\in I$. As a homogeneous system of
$n-\vert I\vert-1$
equations with $n-\vert I\vert$ variables, it must have a non-trivial
solution. Since, by the assumption, $a'_{lr}\neq 0$ for every $r\not\in I$,
we see that this solution gives rise to an invertible diagonal matrix $d$.
Theorem~3 for skew fields is proved.

\smallskip

$2^o$.~It is clear that the assertion of Theorem~3 holds true for
simple Artinian rings, i.e. full matrix rings over skew fields.
And the generalisation of this result to semisimple Artinian rings
is also immediate.

\smallskip

$3^o$.~Let now $R$ be a general Artinian ring. Denote by $G_J$ the principal
congruence subgroup of $G$ modulo $J$, i.e. the subgroup consisting of the
matrices $b\in G$ such that $b_{ij}\equiv e_{ij}(\mod\,J)$. Since the
factor-ring $R/J$ is semisimple Artinian, for every $a\in G$ there
exists a matrix from $\langle a,D\rangle$ which belongs to $G_J$.
So we can suppose that $a\in G_J$, and, further, that $a_{ij}=0$ for
every $i<l$ and $j\neq i$. If $d\in D$, then

$$(ada)_{lj}=\sum\limits_{r=1}^na_{lr}d_ra_{rj}=
a_{ll}d_la_{lj}+a_{lj}d_ja_{jj}+\sum\limits_{r\neq l,j}a_{lr}d_ra_{rj}$$
for every $j\neq l$. Since $a\in G_J$, all its diagonal entries
are invertible. Put $d_l=a_{ll}^{-1},\ d_j=-a_{jj}^{-1}$ for every
$j\neq l$. Then $(ada)_{lj}=\sum\limits_{r\neq l,j}a_{lr}d_ra_{rj}\in
J^2$. Since the radical of the Artinian ring $R$ is nilpotent,
we can argue by induction on the nilpotency exponent in order to
obtain a matrix from
$\langle a,D\rangle$ whose rows from first to $l^{\th}$ coincide with
the corresponding rows of the identity matrix. Now proceed by induction
on $l$. Theorem~3 is proved in full generality.\qed

\sectionheadline{Closing remarks}

$1^o$.~It is known [7] that a cancellative semigroup $(S,\circ)$ is a group
if and only if there exists a function $f:S\to S$ such that the range
of the function $s\mapsto s\circ f(s)$ is finite. Unfortunately, in our
case this elegant result is not efficient.

\smallskip

$2^o$.~Note that the description of intermediate subgroups given in
{\bf Introduction} was obtained for general semilocal rings, not only
for Artinian ones. Therefore it is an open question if Theorem~3 is true
for all semilocal rings (positive answer is plausible since semilocal
rings also contain a lot of invertible elements).

\smallskip

$3^o$.~If $R=k$ is a field (which is sufficiently large), then
from Theorems~1~and~3 (and some results of [2]) it follows that all
intermediate sub(semi)\-groups of $G$ containing $D$ are {\it algebraic},
i.e. they are just sets of $k$--rational points of Zariski closed
sub(semi)groups in $\GL(n,\overline{k})$, where $\overline{k}$ is the algebraic
closure of $k$. The direct proof of this fact is still lacking (even for
subgroups), and it is worth mentioning that each Zariski closed
subsemigroup of an algebraic group is a subgroup, see, e.g. [6].

\smallskip

$4^o$.~If, again, $k$ is a field, then the group $D=\D(n,k)$ is equal to the
group of $k$--rational points of the group $\D(n,\overline{k})$, which is
a maximal split torus in $\GL(n,\overline{k})$. Therefore it is natural to
consider lattices of subsemigroups of $G$ containing groups of $k$--rational
points of other tori in $\GL(n,\overline{k})$, see $5^o, 6^o$ below.
We can also replace here $\GL$ by an algebraic group.
It should be noted that Theorem~3 is not true for the special linear
group $\SL(2,\R)$ and the group $\SD(2,\R)$ consisting of diagonal
matrices with the determinant 1.

\smallskip

$5^o$.~It is easy to verify that for a commutative Artinian ring $R\/$
the same arguments as in the proof of Theorem~3 allow
to establish that each intermediate subsemigroup of $G$
containing the group of diagonal matrices with one fixed entry
equal to $1$ is a subgroup. On the other hand, for the group of
diagonal matrices with two fixed entries equal to $1$ it is not true.

\smallskip

$6^o$.~Let $K/k$ be a finite separable field extension of degree $n$.
Consider the regular embedding of $K$ into the full matrix ring $\M(n,k)$
and denote the image of the multiplicative group $K^*$ of $K$
under this mapping by $\T$. The group $\T$ is equal to the group
of $k$--rational points of a maximal non-split torus in $\GL(n,\overline{k})$.

The complete description of the subgroup lattice $Lat(\T,\GL(n,k))$
has been obtained for certain fields (e.g. local, finite),
but for arbitrary fields only the case $n=2$ has been considered;
see [5] and references in [10]. It turns out that the structure of
the lattice of intermediate subgroups for non-split tori is strikingly
different from the split case, but, nevertheless, for $n=2$
the lattice of intermediate subsemigroups
coincides with the lattice of intermediate subgroups. Indeed, let $g\in G$,
then it follows from Lemma~1~[5] that there exist $t,t'\in\T$
such that $g^{-1}=tgt'$.

\sectionheadline{Acknowledgements}

Discussions concerning this article were initiated during a talk
between the author, Nikolai Vavilov, Boris Lurje and Oleg Izhboldin in
early February~1999. The author is grateful to all participants
of that conversation for their valuable comments.

This research has been carried out in the framework of the
program ``Young Scientists''--1999.

\references

\lastpage
\bye